\documentclass[12pt]{article}
\usepackage{amsmath,amsthm,amssymb}
\usepackage{amssymb}
\usepackage{amsfonts}
\usepackage{graphicx}
\usepackage{hyperref}
\usepackage{bm}
\usepackage[latin1]{inputenc}

\DeclareMathOperator*{\argmin}{arg\,min}

\newtheorem{theorem}{Theorem}[section] 
\newtheorem{lemma}{Lemma}[section]
\newtheorem{cor}{Corollary}[section]
\newtheorem{remark}{Remark}[section]
\newtheorem{definition}{Definition}[section]
\newtheorem{proposition}{Proposition}[section]

\newcommand{\D}{{\rm d}}
\newcommand{\dx}{\, \D x}
\newcommand{\dy}{\, \D y}
\newcommand{\dz}{\, \D z}

\newcommand{\texts}{\textstyle}

\newcommand{\nz}{\mathbb{N}}
\newcommand{\ts}{T}
\newcommand{\gt}{G}

\title{Existence Theory for the EED Inpainting Problem}
\author{Michael Bildhauer \and Marcelo C\'ardenas
  \and Martin Fuchs \and Joachim Weickert}
\date{}

\begin{document}
\parindent0ex
\maketitle

\begin{center}
 {\em Dedicated to Professor Nina Uraltseva on the occasion of her 85th
 birthday.}
\end{center}

\vspace{1mm}

\begin{abstract}
We establish an existence theory for an elliptic boundary value problem 
in image analysis known as edge-enhancing diffusion (EED) inpainting.
The EED inpainting problem aims at restoring missing data in an image 
as steady state of a nonlinear anisotropic diffusion process where
the known data provide Dirichlet boundary conditions.
We prove the existence of a weak solution by applying the Leray-Schauder 
fixed point theorem and show that the set of all possible weak 
solutions is bounded. 
Moreover, we demonstrate that under certain conditions the sequences 
resulting from iterative application of the operator from the 
existence theory contain convergent subsequences.
\end{abstract}



\section{Introduction}
 
Restoring missing data in an EED inpainting problem means that we are given an open set 
$\Omega \subset \mathbb{R}^2$ as the area of an image 
and that the image data are just partially known, i.e.~we consider 
a subset $K \subset \Omega$ and a function $f$: $K \to [0,1]$. Here the values of $f$ represent the 
grey level between pure white and black.\\
 
As a typical example we may imagine that $K$ is a finite union of tiny regions like it is 
illustrated in Figure \ref{fig:uraltseva} below.\\

The complete image is recovered by filling in the 
missing information at $G:=\Omega\setminus K$.\\ 

To this purpose, the boundary value problem
\begin{align}\label{BV_problem1}
\mbox{div}\left(D(\nabla u_\sigma)\nabla u\right)= 0&\quad 
\mbox{on} \quad G \, , \\\label{BV_problem2}
u=f& \quad \mbox{on} \quad \partial K \, , \\\label{BV_problem3}
D\left(\nabla u_\sigma\right)\nabla u\cdot 
\mathcal{N} =0&\quad \mbox{on}\quad 
\partial\Omega\setminus \partial K
\end{align} 
is studied, where $\mathcal{N}$ denotes the outward unit normal to $\Omega$ and where 
for any given function $w\in L^1(G)$ a Gaussian-smoothed version of $w$
with some {\it fixed parameter} $\sigma$ is denoted by $w_\sigma$.\\

In fact, this problem serves as a mathematical model for the anisotropic diffusion 
image inpainting problem \cite{We94e,WW06,GWWB08}.\\

The Dirichlet boundary condition \eqref{BV_problem2} enforces the inpainted result $u$ 
to be coherent with the known data $f$ at the boundary of the known data 
region $K$. The problem is also supplemented with the natural 
Neumann boundary condition \eqref{BV_problem3}.\\

The \textit{diffusion tensor} 
$D$: $\mathbb{R}^2\rightarrow\mathbb{R}^{2\times 2}$  
depending on the gradient of $u_\sigma$ is usually designed 
in order to steer the 
diffusion process in such a way that important geometrical information is 
taken into account.\\

Before going into details, we have to  clarify our notation and assumptions regarding problem 
\eqref{BV_problem1}--\eqref{BV_problem3}.

\paragraph{Notation.}
In the following we always suppose $\Omega \subset \mathbb{R}^2$ to be a bounded Lipschitz domain.
For the definition of the Lebesgue- and Sobolev spaces  $L^{p}(\Omega)$, 
$W^{k,p}(\Omega)$ and its variants we refer to the monograph  \cite{A75}.\\

If $\sigma$ is a given positive parameter, then throughout this paper we denote for any 
function $w\in L^1(G)$ the globally smoothed function with the symbol $w_\sigma$,
\begin{equation}\label{gauss1}
w_\sigma(x) := \int_G k_\sigma (x-z) w(z)\dz  
\quad \mbox{for all}\quad x\in\mathbb{R}^2 \, ,
\end{equation}
where $k_\sigma$ denotes the Gaussian kernel 
\begin{equation}\label{gauss2}
 k_\sigma(z) := \frac{1}{2\pi\sigma^2}\exp
\left(\frac{-|z|^2}{2\sigma^2} \right)\quad \mbox{for all}\quad 
z\in \mathbb{R}^2\, .
\end{equation}

We use the symbol $c$ to denote a generic constant which may change from 
line to line. Whenever we want to stress its dependence on other values 
$\alpha_1, \alpha_2,...$, we write $c(\alpha_1, \alpha_2, ...).$

\paragraph{Main assumptions.} 
\begin{enumerate}
\item We consider a finite number of open disjoint Lipschitz domains $C_i \subset \Omega$, 
$C_i \not= \Omega$, $i=1$,
\dots , $N$, s.t.~${\rm dist} (C_i,C_j)  >0$ for any $i\not= j$, $i$ $j\in \{1,\dots ,N\}$. We then define
\[
C := \bigcup_{i=1}^N C_i \, ,\quad K:= \overline{C}\, ,\quad G:= \Omega \setminus K\, .
\]
For the sake of simplicity we assume throughout this paper that $K \Subset \Omega$,
hence $\partial K \subset \partial G$.
Otherwise we would have to decompose $\partial K = \Gamma_1 \cup \Gamma_2$ with
\[
\Gamma_1 := \partial K \cap \partial G\, , \quad \Gamma_2 = \partial K - \Gamma_1 \, .
\]
\item We suppose that the function $f$ : $C\rightarrow [0,1]$ is of class $W^{1,2}(C)$. Extending this
function to a function of class $W^{1,2}(\Omega)$, we use the same symbol $f$ with a slight abuse of notation. 

\item The diffusion tensor 
$D: \mathbb{R}^2\rightarrow\mathbb{R}^{2\times 2}$ is assumed to be a continuous function
satisfying the symmetry condition  
\[
d_{\alpha,\beta}(p) = 
d_{\beta,\alpha}(p)\quad \mbox{with}\quad D(p) = (d_{\alpha,\beta}
(p))_{1\leq\alpha\leq\beta\leq2}\quad \mbox{for all}\quad  p\in\mathbb{R}^2\, .
\]
Moreover, our hypothesis on non-degenerate ellipticity reads as: there exist positive constants 
$\lambda$, $c_1$, $c_2$ such that  
\begin{equation}\label{uniform_ellipticity}
c_1\left(1+\frac{|p|^2}{\lambda^2}\right)^{-\frac{1}{2}}|q|^2
\leq D(p)q\cdot q\leq c_2|q|^2 \quad
\mbox{for all}\quad p,\, 
q\in \mathbb{R}^2 \, .
\end{equation}
\end{enumerate}

\paragraph{EED inpainting.}
The choice of the ellipticity condition 
\eqref{uniform_ellipticity} covers the inpainting problem 
with \textit{Edge-Enhancing Diffusion (EED)} \cite{We94e,WW06,GWWB08}. \\

As an easy example, $D$ may look like (in fact its representation along the
edge direction $\nabla u^\perp$ 
and across the edge direction, respectively, and the notion "edge direction" 
is motivated by the boundary of level sets of $u$)
\[
D\big(\nabla u_\sigma\big) := 
\left(\begin{array}{cc} 1&0\\[1ex]0&\big(1+|\nabla u_\sigma|^2\big)^{-1/2}\end{array}\right)\,   .
\]

In other words, the diffusion tensor is such that the differential operator
encourages diffusion along edges over diffusion across edges.\\ 

More precisely, $D(\nabla u_\sigma)$ is defined as a $2\times 2$ 
matrix with eigenvectors 
$\nabla u_\sigma^\perp$ and $\nabla u_\sigma$ having as 
corresponding eigenvalues  $1$ and $g(|\nabla u_\sigma|^2),$ where $g$ is
the Charbonnier diffusivity \cite{CBAB94} given by 
$g(s):=\left(1 + s^2 / \lambda^2 \right)^{-\frac{1}{2}}.$\\

Note that the inpainted image $u$ crucially depends on the steering parameter
$\sigma$. Concerning this, a rigorous analysis remains an open problem.\\

However, for the rest of this work we may assume w.l.o.g. that $\lambda=1$ 
except when stated explicitly otherwise.\\

Edge-enhancing anisotropic diffusion was first proposed
as a parabolic evolution process for denoising \cite{We94e},
where it can be regarded as an anisotropic alternative to isotropic
regularisations \cite{CLMC92} of the Perona-Malik filter \cite{PM90}. This 
evolution has been shown to be well-posed in the continuous, space discrete, 
and fully discrete setting \cite{We97}. Later on, the elliptic steady state
equation of EED has been supplemented with Dirichlet data and used for 
inpainting missing regions in matrix-valued images \cite{WW06}.\\ 

Its main application today is inpainting-based lossy image 
compression, where only a sparse, carefully optimised subset of the data 
is stored and the missing data are recovered by EED inpainting \cite{GWWB08}. 
In this context, experiments have shown that EED gives state-of-the-art
results that outperform other partial differential equations in terms of 
reconstruction 
quality \cite{SPME14}.\\ 

Fig.~\ref{fig:uraltseva} depicts an example for 
EED inpainting of sparse data, using a photo of Professor Nina Uraltseva.
In spite of its qualities in practical applications, there is no existence 
theory for EED inpainting so far. Our paper closes this gap.


\begin{figure}[t]

\begin{center}
 \begin{tabular}{ccc}
 \includegraphics[width=0.3\linewidth] {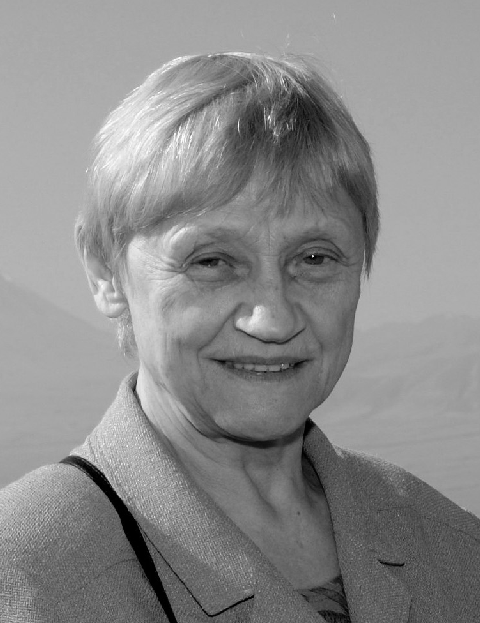} &
 \includegraphics[width=0.3\linewidth] {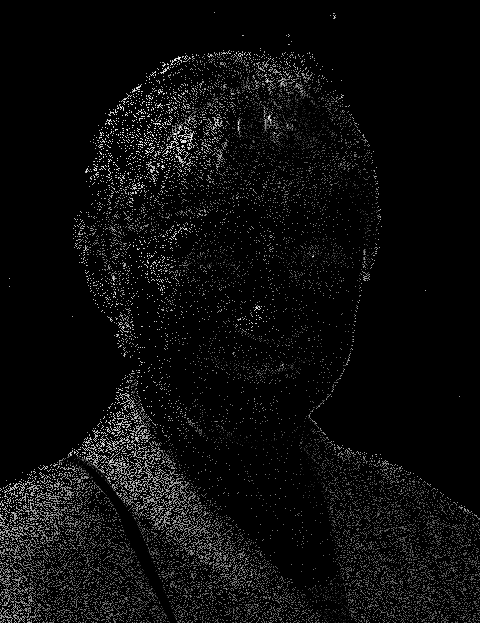} &
 \includegraphics[width=0.3\linewidth] {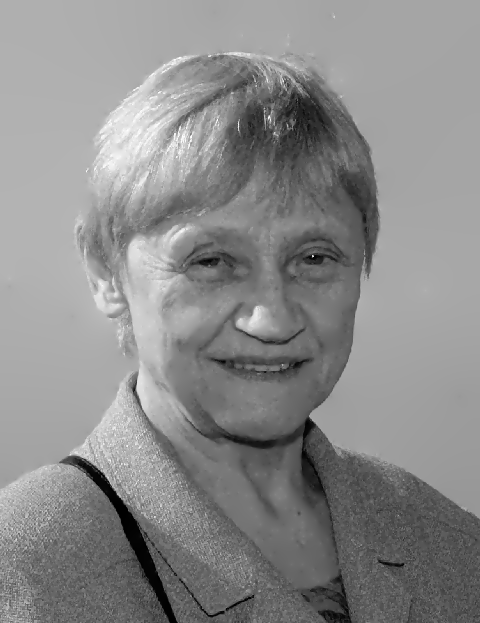}\\
 original image &
 data kept (10 \%) &
 inpainted by EED
 \end{tabular}
\end{center}

\vspace{-5mm}
\caption{\label{fig:uraltseva}
 Illustration of EED inpainting.
 {\bf Left:} Original image courtesy of Professor Uraltseva,
 $\Omega = (0,480) \times (0,623)$.
 {\bf Middle:} Selecting 10 \% of the data with the probabilistic
 sparsification strategy from \cite{MHWT12}. These data locations
 specify the set $K$.
 {\bf Right:} Reconstruction from the sparsified data using a finite
 difference scheme for the EED inpainting model with $\lambda=1$ and
 $\sigma=0.8$.
}

\end{figure}


\paragraph{Main goal and summary.} 
In the present work, we show the existence of a weak solution  
to problem \eqref{BV_problem1}--\eqref{BV_problem3} and  
present some properties of the set of possibly multiple weak solutions.
Moreover, the canonical iteration procedure is analyzed by the way
proving $W^{1,1}$-bounds in different settings.\\

Our paper is organized as follows: Section 2 covers the existence of 
a solution to \eqref{BV_problem1}--\eqref{BV_problem3}. 
We define an operator $T$: $W^{1,1}(G) \to W^{1,1}(G)$ via an associated
variational problem and the Leray-Schauder fixed point theorem 
\cite{GT83,LU68} gives a fixed point of $T$ and in turn the
desired weak solution.\\ 

The analysis of the fixed point set outlined in Section 3 is motivated by the fact, that
no results on the uniqueness of fixed points are available. Here we characterize the fixed points by
another variational formulation.\\

The last section follows the idea, that suitable iterations of the operator $T$ may serve as
a good approximation for a fixed point of $T$.

For this kind of iteration, we shortly establish $W^{1,1}$-bounds in the case of sufficiently large $\sigma$ which 
in general fail to be true. Thus, uniform a priori estimates are based on a suitable smallness assumption.

Given these $W^{1,1}$-bounds, we are able to show the existence of fixed points 
w.r.t.~an iterated operator $T^n$. We emphasize that these fixed points are found in 
a ball whose diameter just depends on the data and is independent of $\sigma$.

Concluding remarks deal with an interpretation of our error estimates.

\section{Existence theory}

The goal of this section is to show that under the above conditions 
problem \eqref{BV_problem1}--\eqref{BV_problem3} has at least one weak 
solution. Our strategy is to apply 
the Leray-Schauder fixed point theorem:\\ 

First we define the weak solutions of the problem as the fixed 
points of an appropriate operator and then show that the operator has  
at least one fixed point \cite{GT83,LU68}. Note that $\sigma$ is fixed throughout this section.\\

Let us begin by introducing the relevant operator. For this purpose 
we need the following two ingredients: the class
\begin{equation}\label{main_space}
\mathcal{C}:=\left\{v \in W^{1,2}(G): \; v=f \; \mbox{on}\; 
\partial K\; \mbox{in the trace sense}\right\} 
\end{equation} 
and the family of functionals $J_w$: $W^{1,2}(G)\rightarrow \mathbb{R}$ 
defined as
$$J_w(v) := \int_G D(\nabla w_\sigma) \nabla v\cdot\nabla v \dx$$
for any $w\in W^{1,1}(G)$ and any $v\in W^{1,2}(G).$

\begin{definition}\label{def_T}
The operator $T: W^{1,1}(G) \rightarrow W^{1,1}(G),$ 
is defined for any $w\in W^{1,1}(G)$ as the argumentum minimi (w.r.t.~the class $\mathcal{C}$)
\begin{equation}\label{functional_problem}
 T(w):=\argmin\limits_{v\in\mathcal{C}} J_w(v) \, .
\end{equation}
For the sake of notational simplicity, the dependence on $\sigma$ is not highlighted unless
we need a careful analysis proving uniform bounds in the last section.
\end{definition}
This definition is well posed on account of Proposition \ref{prop_linear_problem}, whose proof follows
from the ellipticity of the diffusion tensor in a well known manner.
We just recall the basic ingredients concerning Sobolev functions which in addition can be used
as precise references throughout our whole exposition.  The continuity of the trace operator follows, e.g., 
from  Theorem 3.4.5 of \cite{M66}, the variant of Poincar\'e's inequality from Theorem 3.6.4 of \cite{M66}

\begin{proposition}\label{functionspace}
Let $G\subset\mathbb{R}^N$ be a bounded open Lipschitz domain.
\begin{enumerate}
\item \textbf{(Continuity of the trace operator)} \\
For any $1\leq p<\infty$ there exists a bounded 
linear operator  
$B:W^{1,p}(G) \rightarrow L^p(\partial G)$ (the trace operator)
such that $B(u) = u|_{\partial G}$ 
whenever $u\in W^{1,p}(G)\bigcap C(\bar{G}).$ 
In particular, there exists a constant 
$c,$ such that
\begin{equation}\label{trace_ineq}
\|Bu\|_{L^p(\partial G)}\leq c\|u\|_{W^{1,p}(G)},
\end{equation}
for any $u\in W^{1,p}(G).$

\item \textbf{(Poincar\'e inequality)} \\
Let  $\Gamma_1\subset\partial G$ with $\mathcal{H}^{N-1}(\Gamma_1)>0$. 
If $u\in W^{1,p}(G)$ with $1\leq p< \infty$ is such that 
$u|_{\Gamma_1}=0$ in the trace sense, then
\begin{equation}\label{poincare}
\|u\|_{W^{1,p}(G)} \leq c\|\nabla u\|_{L^{p}(G)}
\end{equation}
for a positive constant $c$ which does not depend on $u.$
\end{enumerate}
\end{proposition}

\begin{proposition}\label{prop_linear_problem}
The set $\mathcal{C}$ is a convex and weakly 
closed subset of $W^{1,2}(G).$ Moreover, for any fixed 
$w\in W^{1,1}(G),$ the functional $J_w$ is weakly lower semicontinuous over 
$\mathcal{C}$ and has a unique minimiser in $\mathcal{C}$. 
The unique minimiser $T(w)$ satisfies for all $v\in \mathcal{C}$
\begin{equation}\label{ineq_funcdef}
\int_G {D}(\nabla w_\sigma)\nabla T(w)\cdot\nabla T(w) 
d{x}\leq 
\int_G {D}(\nabla w_\sigma)\nabla v\cdot\nabla v 
\dx
\end{equation}
and for all $\phi\in W^{1,2}(G),$ such that $\phi|_{\partial K} = 0$ the Euler-Lagrange equation
\begin{equation}\label{euler_equation}
\int_G {D}(\nabla w_\sigma)\nabla T(w)\cdot \nabla\phi \dx  = 0\, .
\end{equation}
\end{proposition}

\paragraph{Weak solution.} Any fixed point of the operator $T$ is a weak solution of 
\eqref{BV_problem1}--\eqref{BV_problem3} in the following sense:
\begin{enumerate}
\item By definition, any fixed point $u$ of $T$ satisfies the Euler-Lagrange equation 
\eqref{euler_equation}, hence $u$ is a weak solution of \eqref{BV_problem1} in $G$.

\item Any fixed point is of class $\mathcal{C}$, thus \eqref{BV_problem2} holds.

\item Suppose that $x_0 \in \partial\Omega$ and that there exists an open neighborhood
$U\subset \mathbb{R}^2$ s.t.~$u$ is globally of class $W^{2,2}$ on $U \cap \Omega$. Then we have
\eqref{BV_problem3} a.e.~on $U \cap \partial \Omega$.
\end{enumerate}

We are now ready to state our main result:
\begin{theorem}\label{theo_existence}
\textbf{(Existence of a weak solution)} 
Suppose that we have our general hypotheses as stated above. Then the operator $T$ has at least 
one fixed point, i.e.~the boundary value problem \eqref{BV_problem1}--\eqref{BV_problem3} is
weakly solvable.
\end{theorem} 

The proof of Theorem \ref{theo_existence} 
is divided into different steps which we present as 
intermediate lemmas.
\begin{lemma}\label{lemma1} 
Let $w\in W^{1,1}(G),$ then
\begin{equation}\label{T_estimate}
\|T(w)\| _{W^{1,1}(G)}\leq c\left(1 + \|w\|^\frac{1}{2}_{L^{1}(G)}\right)
\end{equation}
for some positive constant $c$ which is independent of $w.$
\end{lemma}

\emph{Proof.} First we estimate $\nabla w_\sigma$ in 
terms of $w$. It holds for any ${x}\in G$ that
\begin{eqnarray}\label{estimate_wsigma1}
 \partial_\alpha w_\sigma({x}) &=&
 \int_G \partial_\alpha k_\sigma (x-y) w(y) \dy
\end{eqnarray}
for $\alpha =1$, $2$. With \eqref{estimate_wsigma1} we find $c = c(G,\sigma)$ s.t.
\begin{equation}\label{estimate_wsigma2}
\|\nabla w_\sigma\|_{L^{\infty}(G)} \leq c\|w\|_{L^{1}(G)} 
\end{equation}
and using the same argument, \eqref{estimate_wsigma2}
is established for all higher order derivatives.\\

Next we discuss the size of $\|T(w)\|_{W^{1,1}(G)}.$ Fixing $w\in W^{1,1}(G)$, 
we have by H\"older's inequality that
\begin{eqnarray}\label{estimate_wsigma3}
\lefteqn{\int_G |\nabla T(w)|\dx}\\ 
  &=& \int_{G}(1+|\nabla w_\sigma|^2)^{-\frac{1}{4}} |\nabla T(w)| (1+|\nabla w_\sigma|^2)^{\frac{1}{4}}\dx \nonumber\\
 &\leq& \left(\int_G(1+|\nabla w_\sigma|^2)^{-\frac{1}{2}}|\nabla T(w)|^2\dx \right)^{\frac{1}{2}}
 \left( \int_{G}(1+|\nabla w_\sigma|^2)^{\frac{1}{2}}\dx\right)^{\frac{1}{2}}\nonumber\\
 &\leq& c\left(\int_G \left<{D}(\nabla w_\sigma)\nabla 
  T(w), \nabla T(w)\right>\dx \right)^{\frac{1}{2}}
  \left(\int_G(1+|\nabla w_\sigma|^2)^{\frac{1}{2}}\dx\right)^{\frac{1}{2}}\, ,\nonumber
\end{eqnarray}
where for the last inequality we applied the ellipticity  
condition \eqref{uniform_ellipticity}. Moreover, choosing $f$ as 
comparison function in \eqref{ineq_funcdef} 
and applying \eqref{uniform_ellipticity} once again, we obtain
\[
A:= \int_G {D}(\nabla w_\sigma)\nabla 
T(w)\cdot \nabla T(w)\dx \leq c\int_G|\nabla f|^2\dx\, ,
\]
whereas on account of \eqref{estimate_wsigma2},
\[
B:= \int_G(1+|\nabla w_\sigma|^2)^\frac{1}{2}\dx 
\leq c\left(1+\|w\|_{L^{1}(G)}\right)\, .
\]
Applying these estimates of $A$ and $B$ to \eqref{estimate_wsigma3}, 
we obtain 
\begin{equation}\label{estimate_nablaT1}
\|\nabla T(w)\|_{L^1(G)}\leq cA^{\frac{1}{2}}B^{\frac{1}{2}} 
\leq c\left(1+\|w\|^{\frac{1}{2}}_{L^{1}(G)}\right) \, .
\end{equation}
Since $T(w)$ is in the class $\mathcal{C}$, thus coinciding with $f$ over 
$\partial K\subset\partial G,$ we may apply the Poincar\'e 
inequality \eqref{poincare} to the difference 
$T(w) - f$ and obtain
\begin{eqnarray*}
\|T(w)\|_{L^1(G)} &\leq& \|f\|_{L^1(G)} + \|T(w) - f\|_{L^1(G)}\\
&\leq&\|f\|_{L^1(G)} + c\|\nabla(T(w)-f)\|_{L^1(G)} \\
&\leq& c\|f\|_{W^{1,1}(G)} + c\|\nabla T (w)\|_{L^1(G)}.
\end{eqnarray*}
Combining this estimate with inequality \eqref{estimate_nablaT1} 
we end up with \eqref{T_estimate}
for a constant $c$ depending on the data $f$, $G$ and $\sigma$ but not depending on 
the function $w\in W^{1,1}(G)$ under consideration.\hfill$\Box$\\

\begin{lemma}\label{lemma2}
The mapping $T$: $W^{1,1}(G)\rightarrow W^{1,1}(G)$ is continuous.
\end{lemma}

\emph{Proof.} Let $w\in W^{1,1}(G)$ and 
$(w_n)_{n\in\mathbb{N}}\subset W^{1,1}(G)$ be such that 
\begin{equation}\label{eq_cont1}
\|w_n - w\|_{W^{1,1}(G)}\rightarrow 0 \quad \mbox{as} 
\quad n\rightarrow\infty \, .
\end{equation}
Moreover, let $u_n =T(w_n)$ and $u=T(w)$. We have to show that 
\begin{equation}\label{eq_cont2}
\|u_n - u\|_{W^{1,1}(G)}\rightarrow 0 \quad 
\mbox{as} \quad n\rightarrow \infty \, .
\end{equation}

In order to prove \eqref{eq_cont2}, we begin with an observation concerning the mollifications:
recalling \eqref{estimate_wsigma2} and \eqref{eq_cont1} we have
\begin{equation}\label{eq_cont3}
\nabla w^n_\sigma\rightarrow\nabla w_\sigma \quad \mbox{uniformly on } G
\quad 
\mbox{as} \quad n\rightarrow \infty \, ,
\end{equation}
where we use the symbol $w^n_\sigma$ instead of $(w_n)_\sigma$.
The ellipticity condition \eqref{uniform_ellipticity}, the continuity of the diffusion tensor 
and \eqref{eq_cont3} imply for all ${x}\in G, {q}\in\mathbb{R}^2$ and for all $n\in \mathbb{N}$
\begin{equation}\label{eq_cont4}
\mu_1|{q}|^2\leq {D}\big(\nabla w^n_\sigma({x})\big){q}\cdot 
{q}\leq\mu_2|{q}|^2
\end{equation}
with constants $\mu_1$, $\mu_2>0$. 

We again turn our attention to the sequence $u_n$ and combine the minimality of $u_n = T(w_n)$
(recall \eqref{ineq_funcdef}, use $f$ as admissible comparison function) 
with inequality \eqref{eq_cont4} to obtain
\begin{equation}\label{eq_cont4.1}
\mu_1\int_G|\nabla u_n|^2\dx\leq
\int_G {D}(\nabla w^n_\sigma)\nabla u_n\cdot 
\nabla u_n\dx\leq\int_G{D}(\nabla w^n_\sigma)\nabla f\cdot 
\nabla f\dx\, .
\end{equation}

Moreover, since $u_n\in\mathcal{C},$ we may apply Poincar\'e's 
inequality \eqref{poincare} in order to get the estimate
\begin{equation}\label{eq_cont4.2}
\|u_n\|_{L^2(G)}\leq\|u_n-f\|_{L^2(G)} + \|f\|_{L^2(G)}
\leq c\Big(1 + \|\nabla u_n\|_{L^2(G)}\Big)\, 
\end{equation}
again with a constant just depending on the data. Inequalities  
\eqref{eq_cont4.1} and \eqref{eq_cont4.2} are put together with the following result 
of the first step:
\begin{equation}\label{eq_cont5}
\sup\limits_{n\in\mathbb{N}}{\|u_n\|_{W^{1,2}(G)}} < \infty\, .
\end{equation}
Then, by the Banach-Alaoglu theorem, there exists a subsequence  
$(u_{n_k})_{k\in\mathbb{N}}$ of $(u_n)_{n\in\mathbb{N}}$ s.t.~for some $\tilde{u}\in W^{1,2}(G)$
\begin{equation}\label{eq_cont6}
u_{n_k}\rightharpoondown  \tilde{u} \quad \mbox{in} \quad W^{1,2}(G)
\quad 
\mbox{as} \quad k\rightarrow \infty \, .
\end{equation}
In the next step we show that in fact $\tilde{u}=u$ which immediately 
gives the convergence of the sequence as a whole:
\begin{equation}\label{eq_cont9}
u_n\rightharpoondown u\quad \mbox{in}\quad W^{1,2}(G)
\quad 
\mbox{as} \quad n\rightarrow \infty \, .
\end{equation}
In fact, since $\mathcal{C}$ is a closed convex subset of $W^{1,2}(G)$, Mazur's lemma 
implies that it is also weakly closed, hence $\tilde{u}\in\mathcal{C}$. 

On the other hand, we have \eqref{euler_equation} w.r.t.~the minimisers $u_{n_k} =T(w_{n_k})$:  it holds that
\begin{equation}\label{eq_cont7}
\int_G {D}(\nabla w^{n_k}_\sigma)\nabla u_{n_k}\cdot \nabla\phi \dx =0
\end{equation}
for any $\phi\in W^{1,2}(G)$ satisfying $\phi|_{\partial K}=0$.
Hence, applying \eqref{eq_cont3} and \eqref{eq_cont6} in 
equation \eqref{eq_cont7} we arrive at the limit (w.r.t.~$n_k$) equation
\begin{equation}\label{eq_cont8}
\int_G {D}(\nabla w_\sigma)\nabla\tilde{u}\cdot\nabla{\phi} \dx =0\quad\mbox{for any $\phi$ as above.}
\end{equation}
By definition $u$ is a solution of \eqref{eq_cont8} as well which, by uniqueness, is just possible provided that
$\tilde{u} = u$ and we have established our claim \eqref{eq_cont9}. 

It remains to improve 
\eqref{eq_cont9} towards 
\begin{equation}\label{eq_cont10}
\nabla u_n\rightarrow \nabla u\quad \mbox{in}\quad L^2(G,\mathbb{R}^2)
\quad 
\mbox{as} \quad n\rightarrow \infty\, ,
\end{equation}
which gives \eqref{eq_cont2} as a byproduct. We have the Euler-Lagrange equations 
\eqref{euler_equation} both for $u=T(w)$ and for $u_n =T(w_n)$ with $\phi := u_n-u$ as admissible test function
\[
\int_G {D}(\nabla w_\sigma) \nabla u\cdot(\nabla u_n-\nabla u)\dx =
\int_G {D}(\nabla w^n_\sigma) \nabla u_n \cdot(\nabla u_n-\nabla u)\dx = 0 \, , 
\]
which, together with the ellipticity \eqref{uniform_ellipticity}, shows for some $c>0$
\begin{eqnarray*}
c\int_G |\nabla u_n - \nabla u|^2 \dx &\leq& \int_G 
{D}(\nabla w_\sigma)
(\nabla u_n - \nabla u)\cdot(\nabla u_n-\nabla u)\dx\\
&=& \int_G({D}(\nabla w_\sigma) - {D}(\nabla w^n_\sigma))\nabla u_n
\cdot (\nabla u_n - \nabla u)\dx \,.
\end{eqnarray*}
We finally observe that ${D}(\nabla w^n_\sigma)\rightarrow {D}(\nabla w_\sigma)$ 
uniformly on $G$ (by the continuity of ${D}$ and by \eqref{eq_cont3}) which,
together with the uniform bound \eqref{eq_cont5}, gives \eqref{eq_cont10} and 
the proof is completed. \hfill$\Box$\\

The last auxiliary lemma reads as:
\begin{lemma}\label{lemma3}
The continuous mapping $T$: $W^{1,1}(G)\rightarrow W^{1,1}(G)$ is compact.
\end{lemma}
\emph{Proof.} Suppose that we are given a  sequence $(w_n)_{n\in\mathbb{N}}\subset W^{1,1}(G)$ s.t.
\begin{equation}\label{eq_compact1}
\sup\limits_{n\in\mathbb{N}}{\|w_n\|_{W^{1,1}(G)}} < \infty \, .
\end{equation}
Letting $u_n:=T(w_n)$ we then have to extract a convergent subsequence $(u_{n_k})_{k\in\mathbb{N}}$,
\begin{equation}\label{eq_compact2}
u_{n_k}\rightarrow\ : \hat{u}\quad \mbox{in}\quad 
W^{1,1}(G)\quad \mbox{as} \quad k\rightarrow \infty \, ,
\end{equation}
where we do not claim that $\hat{u}$ is a solution of some limit equation.\\

An account of our hypothesis \eqref{eq_compact1} we can exactly follow the lines outlined 
in the proof of Lemma \ref{lemma2} and reproduce \eqref{eq_cont5} in the situation at hand. 
Passing to a subsequence we may suppose
\begin{equation}\label{eq_compact1.1}
u_{n_k} \rightharpoondown: \hat{u}\quad\mbox{in $W^{1,2}(G)$}\, ,
\quad u_{n_k} \to \hat{u} \quad\mbox{in $L^{s}(G)$ for any $1< s <\infty$}\, .
\end{equation}
We finally establsh $L^2$-convergence of the gradients, which in particular 
proves the claim \eqref{eq_compact2}.

In fact, observe that on account of \eqref{estimate_wsigma2} 
and \eqref{eq_compact1} the functions 
\[
{A}_{n_k}:\,\, \bar{G}\rightarrow\mathbb{R}^{2\times 2}\, , \quad
{A}_{n_k}({x}):= {D}(\nabla w^{n_k}_\sigma({x}))\, , 
\]
are bounded and equicontinuous on $\bar{G}$. Thus by Arzela's theorem 
there exist a further subsequence, with a slight abuse of notation still denoted by 
$({A}_{n_k})_{k\in\mathbb{N}}$, and a function  
${A}$: $\bar{G}\rightarrow\mathbb{R}^{2\times 2}$ s.t.~as $k\to \infty$
\[
A_{n_k} \rightrightarrows A\quad\mbox{uniformly in $\bar{G}$.}
\]
We write
\begin{eqnarray}\label{eq_compact4}
\lefteqn{\int_{G}{A}({x})(\nabla u_{n_k} - \nabla\hat{u})\cdot(\nabla u_{n_k} - \nabla\hat{u}) \dx}\\ \nonumber 
&=&\int_{G} {A}_{n_k}({x})(\nabla u_{n_k} - \nabla\hat{u})\cdot 
(\nabla u_{n_k} - \nabla\hat{u})\dx \\ \nonumber
&& + \int_{G}({A}({x})-{A}_{n_k}({x}))
(\nabla u_{n_k} - \nabla\hat{u})\cdot (\nabla u_{n_k} - \nabla\hat{u})\dx \\
&=:& \alpha_k + \beta_k\, . \nonumber
\end{eqnarray}
The uniform convergence ${A}_{n_k}\rightrightarrows {A}$ implies (recalling \eqref{eq_cont5})
\begin{equation}\label{eq_compact5}
\lim_{k\rightarrow\infty}\beta_k =0
\quad 
\mbox{as} \quad k\rightarrow \infty\, .
\end{equation}
As in the proof of Lemma \ref{lemma2}, the Euler-Lagrange equation
w.r.t.~$u_{n_k}=T(w_{n_k})$ is applied leading to
\begin{eqnarray*}
\alpha_k &=& -\int_{G}{A}_{n_k}({x})\nabla\hat{u}\cdot (\nabla u_{n_k} - \nabla \hat{u})\dx \\
&=& \int_{G}\big({A}(x)-{A}_{n_k}({x})\big)\nabla\hat{u}\cdot (\nabla u_{n_k} - \nabla \hat{u})\dx \\
&& -\int_{G}{A}({x})\nabla\hat{u}\cdot (\nabla u_{n_k} - \nabla \hat{u})\dx \, .
\end{eqnarray*}
Here the second term is handled with the weak convergence of $u_{n_k}$ (see \eqref{eq_compact1.1}),
for the first one we argue in addition with the uniform convergence 
 $A_{n_k}\rightrightarrows A$ to obtain
\begin{equation}\label{eq_compact6}
\lim_{k\rightarrow\infty}\alpha_k = 0 
\quad 
\mbox{as} \quad k\rightarrow \infty
\end{equation}
Summarizing the results, \eqref{eq_compact4}, \eqref{eq_compact5} and 
\eqref{eq_compact6} yield the main equality
\begin{equation}\label{eq_compact6.1}
\lim_{k\rightarrow\infty}\int_{G}{A}({x})(\nabla u_{n_k} - \nabla\tilde{u})
\cdot (\nabla u_{n_k} - \nabla\tilde{u}) \dx = 0\, .
\end{equation}
As a final remark we just mention that the condition
\[
\mu_1|{q}|^2\leq {A}({x}){q}\cdot {q}
\]
for all ${x}\in G$, ${q}\in\mathbb{R}^2$ and for some positive constant  $\mu_1$
is induced by the uniform convergence, which completes the proof. \hfill$\Box$

\paragraph{Proof of Theorem \ref{theo_existence}.} The theorem is proved once we have verified the
assumptions of \cite{GT83}, Theorem 11.3.\\ 

By Lemma \ref{lemma3}, $T$ is compact mapping and it remains
to prove the existence of a positive constant $M$ with the property
\[
\|u\|_{W^{1,1}(G)} < M\quad 
\mbox{for all $u\in W^{1,1}(G)$ and for all $\alpha\in[0,1]$ s.t.~$u = \alpha T(u)$.}
\]

Lemma \ref{lemma1} in particular shows
\[
\|T(u)\|_{W^{1,1}(G)}\leq c\left(1 + \|u\|^{\frac{1}{2}}_{W^{1,1}(G)}\right) \, .
\]
W.l.o.g.~we may assume $\alpha>0$. If $u= \alpha T(u)$, then
\[
\|u\|_{W^{1,1}(G)}\leq \alpha c\left(1 + \|u\|^{\frac{1}{2}}_{W^{1,1}(G)}\right)
\]
and using Young's inequality, we obtain for any $\epsilon>0$
\[
\left(1-\frac{\epsilon\alpha c}{2}\right)
\|u\|_{W^{1,1}(G)}\leq \alpha c\left(1 + \frac{1}{2\epsilon}\right)\, .
\]
By elementary calculations we may choose
\[
M = \alpha c\,  \frac{1+2\epsilon}{\epsilon (2-\epsilon \alpha c)} 
\]
letting in addition $\epsilon =1/c$ and the claim follows. \hfill$\Box$

\section{Analysis of the fixed point set}
In the previous section we proved the existence of at least one 
solution to \eqref{BV_problem1}--\eqref{BV_problem3} using the Leray-Schauder 
fixed point argument. \\

Since no information on the uniqueness of solutions is available, the analysis 
of the fixed point set is of particular interest.
In this section we study the properties of the set of possibly multiple weak solutions.\\

We still fix a smoothing parameter $\sigma >0$ throughout this
section.\\ 

In Proposition \ref{prop_Fbounded} we show
that the set
\[
\mathcal{F} := \left\{ u\in W^{1,1}(G):\; u=T(u)\right\} 
\]
of fixed points of $T$ is bounded.\\ 

Then, in Proposition \ref{prop_Jminseq}, we introduce an appropriate functional $\mathcal{J}$
and characterize the set $\mathcal{F}$ as the set of all weak $W^{1,2}$-limits of
$\mathcal{J}$-minimizing sequences.

\begin{proposition}\label{prop_Fbounded}
The set $\mathcal{F}$ is a bounded subset of $W^{1,2}(G)$.
\end{proposition}

\emph{Proof.} By definition we have $T(W^{1,1}(G))\subset \mathcal C$ (recall \eqref{main_space}).
Thus, Poincar\'e's inequality \eqref{poincare} can be applied to any
function $w\in\mathcal{F}$ and we find a constant
$c=c(\sigma,f,G)$ s.t.
\begin{eqnarray*}
\|w\|_{L^1(G)}&\leq& \|w-f\|_{L^1(G)} + \|f\|_{L^1(G)}\\
&\leq& c\left(1 +\|\nabla w\|_{L^1(G)}\right) \nonumber \, .
\end{eqnarray*}
Lemma \ref{lemma1} yields for $w=T(w)$
\[
\|\nabla w\|_{L^1(G)} = \|\nabla T(w)\|_{L^1(G)} \leq c\left(1 + 
\|\nabla w\|^\frac{1}{2}_{L^1(G)}\right)\, .
\]
Youngs' inequality then immediately gives
\begin{equation}\label{kappa}
\|\nabla w\|_{L^1(G)}\leq c(\sigma,f,G) =: \kappa <\infty\quad\mbox{for all $w\in\mathcal{F}$}\, .
\end{equation}
Note that \eqref{kappa} implies as above
\[
\|w\|_{L^1(G)} \leq c(\kappa) \, .
\]
Finally we may apply \eqref{estimate_wsigma2} once again together 
with the ellipticity condition \eqref{uniform_ellipticity}, which
proves the proposition on account of
\[
\int_G |\nabla w|^2\dx \leq c(\kappa) \int_G D(\nabla w_\sigma(x)) \nabla w \cdot \nabla w\dx
\leq c(\kappa) \|\nabla f\|^2_{L^2(G)} \, ,
\]
where the last inequality follows from the minimality of $w=T(w)$. \hfill$\Box$\\

Next we consider the functional 
\[
\mathcal{J}:\,\, \mathcal{C}\rightarrow\mathbb{R}\, , \quad
\mathcal{J}[u] := \int_G|\nabla u - \nabla T(u)|^2 \dx \
\]
and the set of cluster points
\begin{eqnarray*}
\mathcal{M}&:=& \Big\{ u \in W^{1,2}(G):\\
&&\mbox{$u$ is the weak $W^{1,2}$-limit of a $\mathcal{J}$-minimizing sequence}\Big\} \, .
\end{eqnarray*}

\begin{proposition}\label{prop_Jminseq}
The sets $\mathcal{F}$ and $\mathcal{M}$ coincide, i.e.~$\mathcal{F} = \mathcal{M}$.
\end{proposition}
\emph{Proof.}  The inclusion $\mathcal{F}\subset \mathcal{M}$ 
follows (considering a fixed point as a contant sequence)
from Proposition \ref{prop_Fbounded} which gives $\mathcal{F} \subset W^{1,2}(G)$. 
Moreover,
we observe that the existence of a fixed point (Theorem \ref{theo_existence})
implies
\begin{equation}\label{infJ}
\inf_{\mathcal{C}} \mathcal{J} = 0 \, .
\end{equation}
Now fix any $\mathcal{J}$-minimizing sequence $(u_n)_{n\in\mathbb{N}}
\subset\mathcal{C}$ with weak $W^{1,2}$-cluster point $u$, i.e.~we suppose after
passing to a subsequence (not relabeled)
\[
u_{n} \rightharpoondown : \hat{u}\quad \mbox{in $W^{1,2}(G)$}
\] 
and claim that $\hat{u}$ is a fixed point of $T$, i.e.~$\hat{u}=T(\hat{u})$.
Note that $\mathcal{C}$ is weakly closed, i.e.~$\hat{u}\in \mathcal{C}$.\\

We have 
\begin{equation}\label{bound1}
\sup_{n\in\mathbb{N}}\|u_n\|_{W^{1,2}(G)} <\infty 
\end{equation}
and by the $\mathcal{J}$-minimizing property of the sequence
\begin{equation}\label{eq_minseq3.1}
\sup\limits_{n\in\mathbb{N}}\|T(u_n)\|_{W^{1,2}(G)} <\infty\, .
\end{equation}

From \eqref{bound1} and \eqref{eq_minseq3.1} we a find 
subsequence $(u_{n_k})_{k\in\mathbb{N}}\subset 
(u_n)_{n\in\mathbb{N}}$ s.t.
\begin{equation}\label{eq_minseq6.1}
u_{n_k} \rightharpoondown : \hat{u}\, ,\quad
T(u_{n_k}) \rightharpoondown : \xi \quad\mbox{in $W^{1,2}(G)$}\, .
\end{equation} 
where we also have $\xi \in \mathcal{C}$. Moreover, by  \eqref{estimate_wsigma2}
and Arzela's theorem, we may assume the uniform convergence on $G$ as $k\to \infty$:
\begin{equation}\label{eq_minseq7}
{D}(\nabla u^{n_k}_\sigma) \rightrightarrows {D}(\nabla \hat{u}_\sigma)\, .
\end{equation}

Arguing with lower semicontinuity (see \eqref{eq_minseq6.1}), the uniform convergence
\eqref{eq_minseq7}, the minimality of $T(u_{n_k})$ and once again with \eqref{eq_minseq7}
we have for any $v\in\mathcal{C}$
\begin{eqnarray*}
\int_G {D}(\nabla\hat{u}_\sigma)\nabla\xi\cdot\nabla\xi \dx &\leq& 
\liminf_{k\rightarrow\infty}\int_G {D}(\nabla\hat{u}_\sigma)\nabla T(u_{n_k})\cdot\nabla T(u_{n_k})\dx\\
&\leq& 
\liminf_{k\rightarrow\infty} \int_G {D}(\nabla u^{n_k}_\sigma)\nabla T(u_{n_k})\cdot\nabla T(u_{n_k})\dx\\
&\leq & \liminf_{n\rightarrow\infty}\int_G {D}(\nabla u^{n_k}_\sigma)\nabla v\cdot\nabla v \dx\\
&\leq& \int_G{D}(\nabla\hat{u}_\sigma)\nabla v \cdot\nabla v\dx \, .
\end{eqnarray*}
This minimizing property of $\xi$ implies by the uniqueness of solutions $\xi = T(\hat{u})$.
Finally pointwise a.e.~convergence of $u_{n_k}$ and $T(u_{n_k})$ show
the claim $\hat{u} = T(\hat{u})$.
\hfill$\Box$

\begin{remark}\label{minimalw12}
In fact, any $\mathcal{J}$-minimizing sequence $(u_n)_{n\in\mathbb{N}}$
satisfies
\[
\sup_{n\in \mathbb{N}} \|u_n\|_{W^{1,2}(G)} < \infty \, ,
\]
which guarantees the existence of weakly convergent subsequences.
\end{remark}
\emph{Proof.} Fix a $\mathcal{J}$-minimizing sequence $(u_n)_{n\in\mathbb{N}}$. The inequality
\[
\|\nabla u_n\|_{L^1(G)} \leq \|\nabla u_n - \nabla T(u_n)\|_{L^1(G)} 
+ \|\nabla T(u_n)\|_{L^1(G)}
\]
together with (recall \eqref{infJ})
\[
\|\nabla u_n - \nabla T(u_n)\|_{L^1(G)} \leq 
|G|^{\frac{1}{2}}
\mathcal{J}[u_n]^{\frac{1}{2}}\rightarrow 0
\quad 
\mbox{as $n\rightarrow \infty$}
\]
yields for any $\epsilon >0$
\[
\|\nabla u_n\|_{L^1(G)} \leq  \|\nabla T(u_n)\|_{L^1(G)} + \epsilon
\]
provided that $n$ is sufficiently large.
Exactly as outlined in the proof of Proposition \ref{prop_Fbounded}
we derive
\[
\|\nabla T(u_n)\|_{L^1(G)} \leq  
c\left(1 + \|\nabla u_n\|^\frac{1}{2}_{L^1(G)}\right)\, ,
\]
and with Young's inequality we get
\begin{equation}\label{bound1}
\sup\limits_{n\in\mathbb{N}}\|u_n\|_{W^{1,1}(G)} <\infty \, .
\end{equation}
On account of \eqref{bound1} we now may follow the concluding remarks proving 
Proposition \ref{prop_Fbounded} with the result
\[
\|\nabla T(u_n)\|_{L^2(G)}^2 \leq \int_G D(\nabla u^n_\sigma) \nabla T(u_n)\cdot \nabla T(u_n)\dx
\leq c\, ,
\]
hence
\[
\sup\limits_{n\in\mathbb{N}}\|T(u_n)\|_{W^{1,2}(G)} <\infty\, .
\]
Since the sequence is $J$-minimizing, Remark \ref{minimalw12} is obvious. \hfill$\Box$

\section{Study of iterated sequences}
In this section we are interested in the question, whether the above considerations
provide a rigorous analytical framework for
an iteration of  the operator $T$ as approximation of the EED inpainting problem
under consideration.\\

More precisely, recalling \eqref{functional_problem} we fix 
$u^{(0)}\in W^{1,1}(G)$ s.t.~$u^{(0)} = f$ on $\partial K$ and let
\begin{equation}\label{iter0.1}
u^{(j+1)} := T \big(u^{(j)}\big)\, , \quad j\in \mathbb{N}_0\, .
\end{equation}
Here and in the following the dependence of the operator $T$ on the parameter $\sigma$
is of particular interest. Nevertheless we keep the notation of the previous sections and
do not highlight the dependence on $\sigma$.\\

In Section \ref{iterlarge} it becomes evident, that for arbitrary values of $\sigma$ we even cannot
expect uniform $W^{1,1}$-estimates for an iteration of $T$.
We just shortly sketch that assuming $\sigma$ to be large enough leads to $W^{1,1}$-bounds
without further constraints on the data.\\

In Section \ref{iter} some refined a priori estimates are presented which lead to fixed points
of the iterated operator $T^n$ in \emph{a priori bounded balls} $\mathcal{U} \subset W^{1,1}(G)$.
However, these estimates rely on a suitable smallness condition on the data.\\

We finish this section with some concluding remark on error estimates.

\subsection{Iterated sequences with large $\sigma$}\label{iterlarge}

\begin{proposition}\label{prop_iter1}
Consider $u^{(j)}$ as given in \eqref{iter0.1}. There exists a constant $\rho >0$ (compare \eqref{const_expl})
s.t.~in case $\sigma^4 > \rho$
\[
\sup_{j\in\mathbb{N}_0}\|u^{(j)}\|_{W^{1,1}(G)}  < \infty \, .
\]
\end{proposition}
\emph{Proof.} 
For any $w\in\mathcal{C}$ we have 
$\nabla w_\sigma({x}) = \int_G \nabla_x k_\sigma(x-z)w({z})\dz$, hence
by the definition off the Gaussian kernel
\begin{equation}\label{eq_iter0}
|\nabla w_\sigma({x})|\leq \frac{1}{2\pi \sigma^{4}} \|w\|_{L^1(G)} \, .
\end{equation}
Let us recall Poincar\'e's inequality \eqref{poincare} by writing $c= c_P$, 
\begin{equation}\label{eq_iter1}
\|w\|_{L^1(G)} \leq c_P\left(\|\nabla(w-f)\|_{L^1(G)} +  \|f\|_{L^1(G)}\right)\, ,
\end{equation}
Specifying the constants occuring in \eqref{estimate_wsigma3} and recalling
the uniform ellipticity \eqref{uniform_ellipticity} we obtain
\begin{equation}\label{eq_iter2}
\|\nabla T(w)\|_{L^1(G)} \leq 2 \left(\frac{c_2}{c_1}\right)^{\frac{1}{2}} \|\nabla f\|_{L^{2}(G)} 
\left(1+\|\nabla w_\sigma \|_{L^1(G)}\right) \, .
\end{equation}
Now let
\begin{eqnarray}\label{const_expl}
\rho_1 &:=& 
2 \left(\frac{c_2}{c_1}\right)^{\frac{1}{2}} \|\nabla f\|_{L^2(G)}\Bigg[1+\frac{|G|}{2\pi \sigma^4} c_P
\Big(\|f\|_{L^1(G)}+\|\nabla f\|_{L^1(G)}\Big)\Bigg]\, ,\nonumber \\
\rho_2 &:=& 2 \left(\frac{c_2}{c_1}\right)^{\frac{1}{2}} \|\nabla f\|_{L^2(G)}\frac{|G|}{2\pi}c_P \, ,\nonumber\\
\rho &:=& \max\{\rho_1,\rho_2\}\, ,
\end{eqnarray}
thus \eqref{eq_iter0}, \eqref{eq_iter1} and \eqref{eq_iter2} yield
\begin{equation}\label{eq_iter4}
\|\nabla T(w)\|_{L^1(G)} \leq \rho \Big(1 +  \sigma^{-4}\|\nabla w\|_{L^1(G)}\Big)\, .
\end{equation}
Applying iteratively \eqref{eq_iter4} to $u^{(j)}$ we obtain that
\begin{equation}\label{eq_iter3.1}
\|\nabla u^{(j)}\|_{L^1(G)}\leq \sum_{k=1}^{\infty} \rho^k \sigma^{-4(k-1)} 
+ \rho^j\sigma^{-4 j}\|\nabla u^{(0)}\|_{L^1(G)} \, .
\end{equation}
Hence, if we assume $\sigma^4 > \rho$, then a final application 
of Poincar\'e's inequality shows the proposition. \hfill$\Box$\\

Using the compactness of the operator $T$ and
the same reasoning as in Lemma \ref{lemma2} we immediately obtain as a byproduct

\begin{cor}\label{cor_iter} 
Consider $u^{(j)}$ as given in \eqref{iter0.1} and suppose that  
$\sigma^4 > \rho$ with $\rho$ given in \eqref{const_expl}).
\begin{enumerate}
\item There exists a subsequence $(u^{(j_k)})_{k\in\mathbb{N}} \subset(u^{(j)})_{n\in\mathbb{N}_{0}}$ 
and a function $u\in W^{1,1}(G)$  such that $u^{(j_k)}\rightarrow u\in W^{1,1}(G)$ as $k\rightarrow\infty$.

\item We have 
\[
\sup_{j\in\mathbb{N}_0}\|u^{(j)}\|_{W^{1,2}(G)}  < \infty \, .
\]
\end{enumerate}
\end{cor}

\subsection{Iterated sequences with small $\sigma$}\label{iter}


Throughout this section, $\sigma$ is some arbitrary (small) fixed 
positive real number. In particular, $\sigma$ 
is not bounded from below by the data of the problem.\\

\subsubsection{A priori estimates}\label{iter a priori}

Here we derive the main tool for the analysis of the sequence $\{u^{(j)}\}$,
where we recall the definition \eqref{iter0.1} of this sequence and once more
emphasize the ellipticity condition \eqref{uniform_ellipticity}, where $\lambda >0$
now is handled as a free parameter.\\

Similar to \eqref{estimate_wsigma3} one observes 
for $w\in W^{1,1}(G)$, $w=f$ on $\partial K$,
\begin{eqnarray*}
\lefteqn{\|\nabla \ts(w)\|_{L^1(G)}}\\ 
&\leq & \int_G \Bigg( \frac{(c_1 \lambda)^{\frac{1}{2}}}
{(\lambda^2+|\nabla w_\sigma|^2)^{\frac{1}{4}}} |\nabla \ts(w)|
\cdot  \frac{(\lambda^2+|\nabla w_\delta|^2)^{\frac{1}{4}}}
{(c_1\lambda)^{\frac{1}{2}}}\Bigg)\dx\\[1ex]
&\leq & \Bigg( \int_G \frac{c_1\lambda}
{(\lambda^2+|\nabla w_\sigma|^2)^{\frac{1}{2}}} |\nabla \ts(w)|^2\dx\Bigg)^
{\frac{1}{2}}\Bigg(\int_G \frac{1}{c_1\lambda} \,   
(\lambda^2 +|\nabla w_\sigma|^2)^{\frac{1}{2}}\dx\Bigg)^{\frac{1}{2}}\\[1ex]
&=:& A^{\frac{1}{2}} \, B^{\frac{1}{2}} \, .
\end{eqnarray*}
Condition \eqref{uniform_ellipticity} and the minimality of $\ts (w)$ give as before
\[
A^{\frac{1}{2}} 
\leq c_2^{\frac{1}{2}} \Bigg(\int_G |\nabla f|^2\dx \Bigg)^{\frac{1}{2}} \, .
\]
Using the elementary inequality
\[
(a+b)^{\frac{1}{2}} \leq a^{\frac{1}{2}} + b^{\frac{1}{2}}\quad
\mbox{for all $a$, $b\geq 0$}
\]
we further obtain (now emphasizing the precise constants)
\begin{eqnarray*}
B^{\frac{1}{2}} &=& \frac{1}{(c_1\lambda)^{\frac{1}{2}}}
\Bigg(\int_G  (\lambda^2+|\nabla w_\sigma|^2)^{\frac{1}{2}}\dx 
\Bigg)^{\frac{1}{2}}
\leq \frac{1}{(c_1 \lambda)^{\frac{1}{2}}} \Bigg( \int_G (\lambda + |\nabla w_\sigma|)\dx\Bigg)^{\frac{1}{2}}\\[1ex]
&\leq& \frac{1}{(c_1\lambda)^{\frac{1}{2}}}\Big(\lambda |G| +  \|\nabla w_\sigma\|_{L^1(G)}\Big)^{\frac{1}{2}}
\leq  \frac{1}{(c_1\lambda)^{\frac{1}{2}}}\Big(\lambda^{\frac{1}{2}} 
|G|^{\frac{1}{2}} + \| \nabla w_\sigma\|^{\frac{1}{2}}_{L^1(G)}\Big)\, .
\end{eqnarray*}
Thus, it is shown that for all $w \in W^{1,1}(G)$, $w=f$ on $\partial K$,
\begin{equation}\label{iter 2}
\|\nabla \ts (w)\|_{L^1(G)} \leq 
c_2^{\frac{1}{2}} \Bigg(\int_G |\nabla f|^2\dx \Bigg)^{\frac{1}{2}} 
 \frac{1}{(c_1\lambda)^{\frac{1}{2}}}
 \Big(\lambda^{\frac{1}{2}} |G|^{\frac{1}{2}} + \| \nabla w_\sigma\|^{\frac{1}{2}}_{L^1(G)}\Big)\, .
\end{equation}
Note that, for instance, the constant occurring in \eqref{estimate_wsigma2} strongly depends on $\sigma$,
thus we now proceed in a different manner using an integration by parts
\begin{eqnarray}\label{uni 1}
\| \nabla w_\sigma\|_{L^1(G)} &=& \Bigg|\int_G  \int_G  \nabla_x k_\sigma(x-z) w(z) \dz \dx\Bigg|\nonumber\\[1ex]
&\leq&  \Bigg|\int_G  \int_G  \nabla_z k_\sigma(x-z) w(z) \dz \dx\Bigg|\nonumber\\[1ex]
&\leq& \|(\nabla w)_\sigma\|_{L^1(G)} + \| w\|_{L^1(\partial G)} \, .
\end{eqnarray}
Next, the well known property
\begin{equation}\label{uni 2}
\|(\nabla w)_\sigma\|_{L^1(G)} = \int_G \Big|\int_{G} k_\sigma (x-z) \nabla w(z)\dz\Big| \dx  
\leq \|\nabla w\|_{L^{1}(G)}
\end{equation}
is used. Denoting the constant in the trace inequality \eqref{trace_ineq} by $c_T$, the
second term on the r.h.s.~of \eqref{uni 1} is handled via (recalling Poincar\'e's inequality)
\begin{eqnarray}\label{uni 3}
\|w\|_{L^1(\partial G)} &\leq& c_T \|w\|_{W^{1,1}(G)}\nonumber\\
&\leq & c_T \Big[\|w-f\|_{W^{1,1}(G)} + \|f\|_{W^{1,1}(G)}\Big]\nonumber\\
&\leq & c_T \Big[ c_P \|\nabla w\|_{L^{1}(G)} + (1+c_P) \|f\|_{W^{1,1}(G)}\Big] \, .
\end{eqnarray}
Altogether, \eqref{uni 1}, \eqref{uni 2} and \eqref{uni 3} imply
\begin{equation}\label{uni 4}
\| \nabla w_\sigma\|_{L^1(G)} \leq
(1+c_T c_P) \|\nabla w\|_{L^1(G)} + c_T(1+c_P) \|f\|_{W^{1,1}(G)}\, .
\end{equation}
Thus, (\ref{iter 2}) and \eqref{uni 4} yield
\begin{eqnarray}\label{iter 5}
\|\nabla \ts (w)\|_{L^1(\gt)}
&\leq& c_2^{\frac{1}{2}} \Bigg(\int_G |\nabla f|^2\dx 
\Bigg)^{\frac{1}{2}}\nonumber\\[1ex]
&&\cdot 
 \frac{1}{(c_1\lambda)^{\frac{1}{2}}}\Bigg[\lambda^{\frac{1}{2}} |G|^{\frac{1}{2}} 
+ (1+c_Tc_P)^{\frac{1}{2}} \|\nabla w\|^{\frac{1}{2}}_{L^1(\gt)}\nonumber \\ 
&&\quad+ c_T^{\frac{1}{2}} (1+c_P)^{\frac{1}{2}}   \|f\|_{W^{1,1}(G)}^{\frac{1}{2}}\Bigg]\, .
\end{eqnarray}
Define the constants
\begin{eqnarray}\label{iter 6}
K_1 &:=& \Bigg(\frac{c_2}{c_1}\Bigg)^{\frac{1}{2}} 
\|\nabla f\|_{L^2(G)} \Bigg( |G|^{\frac{1}{2}} 
+ \frac{c_T^{\frac{1}{2}}}{\lambda^{\frac{1}{2}}}
(1+c_P)^{\frac{1}{2}} \|f\|_{W^{1,1}(G)}^{\frac{1}{2}}\Bigg)\, ,\\[1ex]
\label{iter 7}
K_2 &:=& \frac{1}{\lambda^{\frac{1}{2}}} 
\Bigg(\frac{c_2}{c_1}\Bigg)^{\frac{1}{2}}\|\nabla f\|_{L^2(G)} (1+c_Tc_P)^{\frac{1}{2}}
\end{eqnarray}
and reformulate (\ref{iter 5}) using \eqref{iter 6} and \eqref{iter 7} as
\begin{equation}\label{iter 8}
\|\nabla \ts(w)\|_{L^1(\gt)} \leq K_1 + K_2 
\|\nabla w\|_{L^1(\gt)}^{\frac{1}{2}}\,,
\end{equation}
which allows us to prove:

\begin{theorem}\label{iter theo 1}
With the notation of above we have for all $j \geq 2$
\begin{equation}\label{iter 9}
\|\nabla u^{(j)} \|_{L^1(\gt)} \leq \sum_{l=1}^j 
K_1^{\texts \frac{1}{2^{l-1}}} K_2^{\texts \sum_{i=0}^{l-2}\frac{1}{2^{i}}}
+K_2^{\texts \sum_{i=0}^{j-1} \frac{1}{2^i}} 
\|\nabla u^{(0)}\|_{L^1(\gt)}^{\texts \frac{1}{2^j}}\, .
\end{equation}
\end{theorem}

\emph{Proof by induction.} For $j=2$ inequality (\ref{iter 8}) implies
\begin{eqnarray*}
\|\nabla u^{(2)}\|_{L^1(\gt)}&=& \|\nabla \ts(u^{(1)})\|_{L^1(\gt)}\\[1ex]
&\leq & K_1 + K_2 \|\nabla u^{(1)}\|_{L^1(\gt)}^{\frac{1}{2}}\\[1ex]
&\leq& K_1 + K_2 \Big[K_1+ K_2\|\nabla u^{(0)}\|_{L^1(\gt)}^{\frac{1}{2}}\Big]^{\frac{1}{2}}\\[1ex]
&\leq& K_1 + K_1^{\frac{1}{2}} K_2 + K_2^{1+\frac{1}{2}} 
\|\nabla u^{(0)}\|_{L^1(\gt)}^{\frac{1}{4}}\, ,
\end{eqnarray*}
which is (\ref{iter 9}) for $j=2$.
Now suppose that (\ref{iter 9}) holds for some $j\geq 2$. 
As above for $k\in \nz$ and for given $a_i\geq 0$, $i=1$, 
\dots $k$, the elementary estimate
\[
\Bigg(\sum_{i=1}^k a_i\Bigg)^{\frac{1}{2}} \leq \sum_{i=1}^k 
a_i^{\frac{1}{2}}
\]
is used and again with the help of (\ref{iter 8}) we obtain 
\begin{eqnarray*}
\lefteqn{\|\nabla u^{(j+1)}\|_{L^1(\gt)} = \|\nabla \ts(u^{(j)})\|_{L^{1}(\gt)}
\leq K_1+K_2 \|\nabla u^{(j)}\|_{L^1(\gt)}^{\frac{1}{2}}}\\[1ex]
&\leq & K_1 + K_2 \Bigg( \sum_{l=1}^j K_1^{\texts\frac{1}{2^{l-1}}} 
K_2^{\texts\sum_{i=0}^{l-2}\frac{1}{2^{i}}}
+K_2^{\texts\sum_{i=0}^{j-1} \frac{1}{2^i}} \|\nabla u^{(0)}\|_
{L^1(\gt)}^{\texts\frac{1}{2^j}}\Bigg)^{\frac{1}{2}}\\[1ex]
&\leq&  K_1 + K_2 \Bigg(\sum_{l=1}^j K_1^{\texts\frac{1}{2^l}}K_2^
{\texts\sum_{i=0}^{l-2}\frac{1}{2^{i+1}}}
+K_2^{\texts\sum_{i=0}^{j-1} \frac{1}{2^{i+1}}} \|\nabla u^{(0)}\|_
{L^1(\gt)}^{\texts\frac{1}{2^{j+1}}}\Bigg)\\[1ex]
&=& K_1 +  \sum_{\tilde{l}=2}^{j+1} K_1^{\texts\frac{1}{2^{\tilde{l}-1}}}
K_2 K_2^{\texts\sum_{i=0}^{\tilde{l}-3}\frac{1}{2^{i+1}}}
+K_2 K_2^{\texts\sum_{i=0}^{j-1} \frac{1}{2^{i+1}}} \|\nabla 
u^{(0)}\|_{L^1(\gt)}^{\texts\frac{1}{2^{j+1}}}\\[1ex]
&=& K_1 +  \sum_{l=2}^{j+1} K_1^{\texts\frac{1}{2^{l-1}}} 
K_2^{\texts1+\sum_{i=1}^{l-2}\frac{1}{2^{i}}}
+K_2^{\texts 1+\sum_{i=1}^{j} \frac{1}{2^{i}}} \|\nabla 
u^{(0)}\|_{L^1(\gt)}^{\texts\frac{1}{2^{j+1}}}\\[1ex]
&=&  \sum_{l=1}^{j+1} K_1^{\texts \frac{1}{2^{l-1}}} 
K_2^{\texts\sum_{i=0}^{l-2}\frac{1}{2^{i}}}
+K_2^{\texts\sum_{i=0}^{j} \frac{1}{2^{i}}} 
\|\nabla u^{(0)}\|_{L^1(\gt)}^{\texts\frac{1}{2^{j+1}}}\,. \qquad\qquad\Box
\end{eqnarray*}

\begin{cor}\label {iter cor 1}
Let us suppose that in addition to the assumptions 
of Theorem \ref{iter theo 1}
we have $K_1 < 1$, where $K_1$ is the constant defined in (\ref{iter 6}). \\

Moreover, fix $\rho >0$ and suppose that $\mathcal{W}$ is a subset of 
$W^{1,1}(G)$ 
such that for all $v\in \mathcal{W}$ we have $v =f$ on $\partial K$ and
\[
\|v\|_{W^{1,1}(G)} \leq \rho \, .
\]
\begin{enumerate}
\item Then there is a uniform constant $c=c(\rho)$ such that for all 
$v^{(0)} \in \mathcal {W}$,
$v = v^{(n)} = T\big(v^{(n-1)}\big)$, $n\in \nz$, we have
\[
\|v\|_{W^{1,1}(\gt)} \leq c\, .
\]
\item For a universal constant $c$ not depending on $\rho$ we have i) 
for any $n \geq n_0 = n_0(\rho)$.  
\end{enumerate}
\end{cor}

\emph{Proof.} 
Once again recalling \eqref{poincare}, the proof follows from Theorem 
\ref{iter theo 1}.\hspace*{\fill}$\Box$\\

\begin{remark}\label{iter rem 1}
As it should be expected, the smallness of $K_1$ corresponds to
$c_2\approx c_1$ or $\|\nabla f\|$ small or $|G|$ small or  $\lambda$ large.
\end{remark}
\subsubsection{Iterated fixed points on a priori sets}\label{iter fix}
In this subsection we derive the existence of iterated fixed points in 
a priori local bounded subsets of $W^{1,1}$.\\

In order to find a suitable domain $\mathcal{M}$ s.t.~$T$: 
$\mathcal{M} \to \mathcal{M}$, we first let
for any fixed $w\in W^{1,1}(G)$
\begin{eqnarray*}
\mathcal{T}[w] &:=& \bigcup_{n=1}^{\infty} T^n(w)\\[1ex]
&=& \left\{ v \in W^{1,1}(G):\, 
v=T^n(w)\quad\mbox{for some $n\in \nz$}\right\}\, .
\end{eqnarray*}
We then fix some $\mathcal{U}\subset W^{1,1}(G)$,
where for our purposes we may suppose in the following 
that $v=f$ on $\partial K$ for any $v\in \mathcal{U}$.\\

Now let
\[
\mathcal{M} := \mathcal{U}\cup \bigcup_{v \in \mathcal{U}} \mathcal{T}[v] 
=\mathcal{U} \cup T(\mathcal{U}) \cup T\big(T(\mathcal{U})\big) \dots \, .
\]

Clearly, if $v\in \mathcal{M}$ then $v\in \mathcal{U}$ or $v=T^{n}(w)$ 
for some $w\in \mathcal{U}$ and some $n\in \nz$,
hence
\[
T(v) \in \mathcal{T}[v] \cup \mathcal{T}[w] \subset \mathcal{M}\, .
\]

With the additional notation
\begin{eqnarray*}
\mathcal{T}_{n_0}[w] &:=& \bigcup_{k= n_0}^{\infty} T^k(w)\, , 
\quad\mbox{$n_0\in \nz$ fixed}\, ,\\[1ex]
\mathcal{M}_{n_0} &:=& \mathcal{U}\cup \bigcup_{v \in \mathcal{U}} 
\mathcal{T}_{n_0}[v]\, , 
\end{eqnarray*}
the same reasoning as above shows for all $n\in \nz$
\[
T^n:\quad \mathcal{M}_n \to \bigcup_{v\in \mathcal{U}} \mathcal{T}_n[v] 
\subset \mathcal{M}_n\, .
\]

Increasing $\rho$ in Corollary \ref{iter cor 1}, $ii)$, if 
necessary, and choosing $n_0$ sufficiently large,
we now choose
\[
\mathcal{U} = \{v\in W^{1,1}(G): \, 
\|\nabla v\|_{L^{1}(\gt)} \leq \rho\, ,\, \mbox{$v=f$ on $\partial K$}\}
\]
such that for all $n \geq n_0$
\[
\bigcup_{v\in \mathcal{U}}\mathcal{T}_{n}[v] \subset 
\mathcal{U} = \mathcal{M}_{n}\, .
\]
Since $\mathcal{U}$ is a closed convex subset 
of $W^{1,1}(G)$, the above equality shows the same for $\mathcal{M}_{n}$,
and Corollary 11.2 of \cite{GT83} provides for all $n \geq n_0$ at least 
one fixed point of $T^n$ in $\mathcal{M}_n$,
hence in $\mathcal{M}_{n_0}$.

\subsubsection{Error estimates}\label{iter error}
The numerical discussion of iterated fixed points is motivated by 
letting for all $n$, $k\in \nz$
and for a given $u^{(0)}$ as above (recall \eqref{iter0.1}
\[
 u^{(n+k)} = T^n T^k \big(u^{(0)}\big) = T^n \big(u^{(k)}\big) \, .
\]
With the notation
\[
R(n,k) := u^{(n+k)}-u^{(k)}
\]
we immediately obtain the first observation:
\[
T^{n}\big(u^{(k)}\big) = u^{(k)} + R(n,k)\, ,
\]
i.e.~$\|R(n,k)\|$ provides a measure for the failure of $u^{(k)}$ to 
be numerically detected 
as one possible fixed point w.r.t.~$T^n$.\\

Of course it is possible to choose subsequences $\{n_j\}$ and 
$\{k_j\}$ such that as $j\to \infty$
(recall Theorem \ref{iter theo 1})
\[
\|R(n_j,k_j)\|_{W^{1,1}(\gt)} \to 0 \, .
\]

These simple observations lead to the analysis of the operator 
$T^n$. In fact, if $w_n$ is a fixed point
w.r.t.~$T^n$,
\[
w_n = T^n(w_n)\, ,
\]
then we have
\begin{eqnarray}\label{iter 10}
T(w_n) &=& T^{n+1}(w_n) = T^n(w_n) + 
\big(T^{n+1}-T^n\big)(w_n)\nonumber\\[2ex]
&=& w_n + T^n \big(T-{\rm id}\big)(w_n) \, .
\end{eqnarray}
Now a given fixed point $w_n$ of $T^n$ in general is not simultaneously a 
fixed point of $T,$ although the
converse trivially is true.
The error $v_n:= (T-id)(w_n)$ evidently has to enter the left hand side 
of (\ref{iter 10}).\\

Note that we then apply the iterated operator $T^n$ to $v_n$. This 
is crucial for the numerical interpretation: Decreasing $\sigma$ in the 
kernel
appearing in $T=\ts$ means to blow up the error via the 
non-contracting operator $\ts$ in $\ts^n$
which exactly describes the behavior of the examples sketched in 
the introduction.

\section*{Acknowledgement}
The research of M.C. and J.W. has received funding by the European Research 
Council (ERC) under the European Union's Horizon 2020 research and 
innovation programme (grant agreement no. 741215, ERC Advanced Grant INCOVID).


\bibliographystyle{splncs04}
\bibliography{refs}


\vspace*{0.5cm}
\begin{tabular}{ll}
Michael Bildhauer*&bibi@math.uni-sb.de\\
Marcelo C\'ardenas{**}&cardenas@mia.uni-saarland.de\\ 
Martin Fuchs*&fuchs@math.uni-sb.de\\ 
Joachim Weickert{**}&weickert@mia.uni-saarland.de
\end{tabular}

\vspace*{0.5cm}
\small
\begin{tabular}{llrl}
*&Department of Mathematics&\quad **&Math.~Image Analysis Group\\
&Saarland University&&Saarland University\\
&Faculty Math.~and Computer Sci.&&Faculty Math.~and Computer Sci.\\
&P.O.~Box 15 11 50&&Campus E1.7\\
&66041 Saarbr\"ucken, Germany&&66041 Saarbr\"ucken, Germany \\
\end{tabular}


\end{document}